\input amstex
\documentstyle{amsppt}
%
\catcode`@=11
\redefine\output@{%
  \def\break{\penalty-\@M}\let\par\endgraf
  \ifodd\pageno\global\hoffset=105pt\else\global\hoffset=8pt\fi  
  \shipout\vbox{%
    \ifplain@
      \let\makeheadline\relax \let\makefootline\relax
    \else
      \iffirstpage@ \global\firstpage@false
        \let\rightheadline\frheadline
        \let\leftheadline\flheadline
      \else
        \ifrunheads@ 
        \else \let\makeheadline\relax
        \fi
      \fi
    \fi
    \makeheadline \pagebody \makefootline}%
  \advancepageno \ifnum\outputpenalty>-\@MM\else\dosupereject\fi
}
\def\Beta{\mathchar"0\hexnumber@\rmfam 42}
\redefine\mm@{2010} 
\catcode`\@=\active
\nopagenumbers
\chardef\textvolna='176

\chardef\bigalpha='013
\def\negskp{\hskip -2pt}
\def\GCD{\operatorname{GCD}}
\chardef\degree="5E
\def\compos{\,\raise 1pt\hbox{$\sssize\circ$} \,}

\def\blue#1{#1}

\gdef\darkred#1{#1}

\catcode`#=11\def\diez{#}\catcode`#=6
\catcode`&=11\catcode`&=4
\catcode`_=11\def\podcherkivanie{_}\catcode`_=8
\catcode`\^=11\catcode`\^=7
\catcode`~=11\catcode`~=\active
\catcode`\%=11\def\procent{
\def\mycite#1{\cite{\blue{#1}}\immediate\special{ps:
     ShrHPSdict begin /ShrBORDERthickness 0 def}}
\def\myciterange#1#2#3#4{\cite{\blue{#2#3#4}}\immediate\special{ps:
     ShrHPSdict begin /ShrBORDERthickness 0 def}}
\def\mytag#1{%
    \tag#1}
\def\mythetag#1{\thetag{\blue{#1}}\immediate\special{ps:
     ShrHPSdict begin /ShrBORDERthickness 0 def}}
\def\myrefno#1{\no#1}
\def\myhref#1#2{\blue{#2}\immediate\special{ps:
     ShrHPSdict begin /ShrBORDERthickness 0 def}}
\def\myEarXivlink{\myhref{http://arXiv.org}{http:/\negskp/arXiv.org}}
\def\myGeoCities{\myhref{http://www.geocities.com}{GeoCities}}
\def\mytheorem#1{\csname proclaim\endcsname{Theorem #1}}
\def\mytheoremwithtitle#1#2{\csname proclaim\endcsname{Theorem #1#2}}
\def\mythetheorem#1{\blue{#1}\immediate\special{ps:
     ShrHPSdict begin /ShrBORDERthickness 0 def}}
\def\mylemma#1{\csname proclaim\endcsname{Lemma #1}}
\def\mylemmawithtitle#1#2{\csname proclaim\endcsname{Lemma #1#2}}
\def\mythelemma#1{\blue{#1}\immediate\special{ps:
     ShrHPSdict begin /ShrBORDERthickness 0 def}}
\def\mycorollary#1{\csname proclaim\endcsname{Corollary #1}}
\def\mythecorollary#1{\blue{#1}\immediate\special{ps:
     ShrHPSdict begin /ShrBORDERthickness 0 def}}
\def\mydefinition#1{\definition{Definition #1}}
\def\mythedefinition#1{\blue{#1}\immediate\special{ps:
     ShrHPSdict begin /ShrBORDERthickness 0 def}}
\def\myconjecture#1{\csname proclaim\endcsname{Conjecture #1}}
\def\myconjecturewithtitle#1#2{\csname proclaim\endcsname{Conjecture #1#2}}
\def\mytheconjecture#1{\blue{#1}\immediate\special{ps:
     ShrHPSdict begin /ShrBORDERthickness 0 def}}
\def\myproblem#1{\csname proclaim\endcsname{Problem #1}}
\def\myproblemwithtitle#1#2{\csname proclaim\endcsname{Problem #1#2}}
\def\mytheproblem#1{\blue{#1}\immediate\special{ps:
     ShrHPSdict begin /ShrBORDERthickness 0 def}}
\def\mytable#1{Table #1}
\def\mythetable#1{\blue{#1}\immediate\special{ps:
     ShrHPSdict begin /ShrBORDERthickness 0 def}}
\def\myanchortext#1#2{#2}
\def\mytheanchortext#1#2{\blue{#2}\immediate\special{ps:
     ShrHPSdict begin /ShrBORDERthickness 0 def}}
\font\eightcyr=wncyr8
\pagewidth{360pt}
\pageheight{606pt}
\topmatter
\title
Evolution patterns in Collatz problem.
\endtitle
\author
Ruslan Sharipov
\endauthor
\address Bashkir State University, 32 Zaki Validi street, 450074 Ufa, Russia
\endaddress
\email
\myhref{mailto:r-sharipov\@mail.ru}{r-sharipov\@mail.ru}
\endemail
\abstract
The concept of evolution patterns is introduced for Collatz sequences and it is shown 
that any finite evolution pattern is implemented in some particular Collatz sequence. 
\endabstract
\endtopmatter
\loadbold
\TagsOnRight
\document

\head
1. Introduction.
\endhead
     Wikipedia \mycite{1} says that the Collatz conjecture also known as the 
$3\,n+1$ problem was formulated by Lothar Collatz in 1937. However Wikipedia 
gives a lot of different names of this conjecture associated with different 
persons and not only with persons: the Ulam conjecture, Kakutani's problem, the
Thwaites conjecture, Hasse's algorithm, the Syracuse problem.\par
     The statement of the Collatz conjecture is based on the mapping $f\!:\,
\Bbb N\to\Bbb N$ which is defined as follows in the set of positive integers
$\Bbb N$: 
$$
\hskip -2em
f(n)=\cases 3\,n+1 &\text{if $n$ is odd},\\
n/2 &\text{if $n$ is even}.
\endcases
\mytag{1.1}
$$
Starting with an arbitrary number $m\in\Bbb N$, the sequence of numbers $a_i$ is built 
such that $a_1=m$ and $a_{i+1}=f(a_i)$ for all $i\in\Bbb N$. For instance in the case  
of $m=1$, we get the sequence of numbers
$$
\hskip -2em
1,\ 4,\ 2,\ 1,\ 4,\ 2,\ 1,\ \ldots\ ,
\mytag{1.2}
$$
which repeats periodically with the period $T=3$. The cases $m\neq 1$ are described by
the Collatz conjecture.\par
\myconjecturewithtitle{1.1}{ (Collatz)} For any positive integer $m\in\Bbb N$ the
sequence of numbers $a_i$ defined using the mapping \mythetag{1.1} through $a_1=m$ 
and $a_{i+1}=f(a_i)$ for all $i\in\Bbb N$ reaches the number $1$ \rm (\,i\.\,e\. $a_k=1$
for some $k\in\Bbb N$\,\rm ) and then repeats periodically as in \mythetag{1.2}. 
\endproclaim
     There are many research works devoted to the Collatz conjecture~\mytheconjecture{1.1}, 
see \mycite{2} and \mycite{3}. The year of 2021 and early 2022 demonstrated rather high 
publication activity in this field, see \myciterange{4}{4}{--}{33}. We are not going to 
analyze all of these papers hear. Our goal is to show that the evolution of integers in 
Collatz sequences can be arbitrarily complicated in such a way that any predefined finite 
evolution pattern can be implemented in some Collatz sequence.\par
\head
2. Evolution patterns.
\endhead
\mylemma{2.1} Any positive even integer $n$ can be uniquely presented as $n=2^s\cdot m$,
where $m$ and $s$ are two positive integers and $m$ is odd. 
\endproclaim
    The proof consists in iteratively dividing $n$ by $2$ until the odd number $m$ is
reached: $n\to n/2\to n/2^2\to\ldots\to n/2^s=m$. 
\mylemma{2.2} Any positive odd integer $n$ can be uniquely presented as $n=2^s\cdot m-1$,
where $m$ and $s$ are two positive integers and $m$ is odd. 
\endproclaim
    The proof consists in applying Lemma~\mythelemma{2.1} to the even number $n+1$.\par
    Let's consider the Collatz evolution initiated by the odd number $a_1=2^s\cdot m-1$.
Due to \mythetag{1.1} we get
$$
\hskip -2em
a_2=f(a_1)=3\,a_1+1=3\cdot 2^s\cdot m-2.
\mytag{2.1}
$$
The number $a_2$ is even, therefore
$$
\hskip -2em
a_3=f(a_2)=a_2/2=3^1\cdot 2^{s-1}\cdot m-1.
\mytag{2.2}
$$
If $s-1\neq 0$, then $a_3$ is again odd and we can repeat the steps
\mythetag{2.1} and \mythetag{2.2}:
$$
\hskip -2em
\aligned
&a_4=f(a_3)=3\,a_3+1=3^2\cdot 2^{s-1}\cdot m-2,\\
&a_5=f(a_4)=a_4/2=3^2\cdot 2^{s-2}\cdot m-1.
\endaligned
\mytag{2.3}
$$
Now due to \mythetag{2.2} and \mythetag{2.3} it is easy to see that
$$
\hskip -2em
a_{1+2s}=3^s\cdot m-1.
\mytag{2.4}
$$
The number \mythetag{2.4} is even. 
\mydefinition{2.1} The Collatz evolution from an odd number $a_1=2^s\cdot m-1$
to the even number $a_{1+2s}=3^s\cdot m-1$ is called the $s$-evolution.
\enddefinition
    Applying Lemma~\mythelemma{2.1} to the even number \mythetag{2.4}, we get
$$
\hskip -2em
a_{1+2s}=3^s\cdot m-1=2^q\cdot\tilde n,
\mytag{2.5}
$$
where $q\geqslant 1$ and $\tilde n$ is a positive odd number. Applying the Collatz
evolution to the even number \mythetag{2.5}, we find
$$
\hskip -2em
a_{1+2s+q}=a_{1+2s}/2^q=\tilde n.
\mytag{2.6}
$$
\mydefinition{2.2} The Collatz evolution from an even number $a_{1+2s}=2^q\cdot\tilde n$
to the odd number $a_{1+2s+q}=\tilde n$ is called the $q$-evolution.
\enddefinition
     Since the number \mythetag{2.6} is again odd, we can apply the steps from
\mythetag{2.1} to \mythetag{2.6} to it. As a result we see that we have proved the 
following theorem.
\mytheorem{2.1} The Collatz evolution of any odd number is presented by an alternating sequence
of $s$ and $q$ evolutions:
$$
\hskip -2em
s_1,\,q_1,\,s_2,\,q_2,\,\ldots,\,s_r,\,q_r,\,s_{r+1},\,\ldots\ .
\mytag{2.7}
$$
\endproclaim  
\mydefinition{2.3} The infinite sequence of positive integers \mythetag{2.7} is called
the Collatz evolution pattern of a given odd number. 
\enddefinition
    For instance the evolution pattern of the odd number $1$ is given by the following
trivial periodic sequence
$$
\hskip -2em
1,\,1,\,1,\,\ldots,\,1\,\ldots\ .
\mytag{2.8}
$$
The sequence \mythetag{2.8} is easily derived from \mythetag{1.2}. 
\mydefinition{2.4} Any finite initial part of the Collatz evolution pattern \mythetag{2.7} 
is called the finite Collatz evolution pattern of a given odd number:
$$
\hskip -2em
s_1,\,q_1,\,s_2,\,q_2,\,\ldots,\,s_r,\,q_r,\,s_{r+1}.
\mytag{2.9}
$$
\enddefinition
     Our further goal is to prove the following theorem.
\mytheorem{2.2} For any finite sequence of positive integers \mythetag{2.9} there
is an odd number whose finite Collatz evolution pattern coincides with \mythetag{2.9}.
\endproclaim
\head
3. Diophantine equations associated with $sqs$-patterns.
\endhead     
     Assume that $r=1$ in \mythetag{2.9}. Then we have the following finite Collatz 
evolution pattern: $s_1,\,q_1,\,s_2$. The actual evolution associated with this pattern 
looks like
$$
\hskip -2em
\CD
2^{s_1}\cdot m_1-1@>s_1>> 3^{s_1}\cdot m_1-1@>q_1>>2^{s_2}\cdot m_2-1.
\endCD
\mytag{3.1}
$$
Relying on \mythetag{3.1}, we can write the following equality:
$$
\hskip -2em
3^{s_1}\cdot m_1-1=2^{q_1}\cdot(2^{s_2}\cdot m_2-1).
\mytag{3.2}
$$
The equality \mythetag{3.2} can be rewritten as
$$
\hskip -2em
2^{q_1+s_2}\cdot m_2-3^{s_1}\,m_1=2^{q_1}-1.
\mytag{3.3}
$$
Since $s_1$, $q_1$, and $s_2$ are known, we can treat \mythetag{3.3} as
a linear Diophantine equation with respect to the variables $m_1$ and
$m_2$.\par
    The theory of linear Diophantine equations is given in Section 2.1
of Chapter I.2 in the book \mycite{34}. There one can find the following 
theorem.
\mytheorem{3.1} Let $a$, $b$, $c$ be integers, $a$ and $b$ nonzero. Consider 
the linear Diophantine equation
$$
\hskip -2em
a\,x+b\,y=c
\mytag{3.4}
$$
\roster
\rosteritemwd=0pt
\item"1." The equation \mythetag{3.4} is solvable in integers if and only if 
the greatest common divisor $d=\GCD(a,b)$ divides $c$.
\item"2." If $(x,y)=(x_0,y_0)$ is a particular solution to \mythetag{3.4}, then
every integer solution is of the form
$$
\xalignat 2
&\hskip -2em
x=x_0+\frac{b}{d}\cdot t,
&&y=y_0-\frac{a}{d}\cdot t,
\mytag{3.5}
\endxalignat
$$
where $t$ is an integer. 
\item"3." If $c=GCD(a,b)$ and if $|a|$ or $|b|$ is different from $1$, then a particular 
solution $(x,y)=(x_0,y_0)$ in \mythetag{3.5} can be found such that $|x_0|<|b|$ and
$|y_0|<|a|$.
\endroster
\endproclaim
    Comparing \mythetag{3.3} with \mythetag{3.4}, we find that in our particular case 
$$
\xalignat 3
&\hskip -2em
a=2^{q_1+s_2},
&&b=-3^{s_1},
&&c=2^{q_1}-1,
\mytag{3.6}
\endxalignat
$$
i\.\,e\. the equation \mythetag{3.3} is written as
$$
\hskip -2em
2^{q_1+s_2}\cdot x-3^{s_1}y=2^{q_1}-1.
\mytag{3.7}
$$
From \mythetag{3.6} we derive
$$
\hskip -2em
d=\GCD(a,b)=1.
\mytag{3.8}
$$
Applying Item~1 from Theorem~\mythetheorem{3.1} to \mythetag{3.8}, we derive 
the following theorem.
\mytheorem{3.2} For any sequence of positive integers \mythetag{2.9} the Diophantine 
equation \mythetag{3.7} is solvable. 
\endproclaim
     Then applying Item~2 from Theorem~\mythetheorem{3.1} to the equation \mythetag{3.7}, 
we find that the formulas \mythetag{3.5} are written as
$$
\xalignat 2
&\hskip -2em
x=x_0-3^{s_1}\,t,
&&y=y_0-2^{q_1+s_2}\,t,
\mytag{3.9}
\endxalignat
$$
The second equality \mythetag{3.9} means that we can choose a unique particular solution 
$(x_0,y_0)$ of the equation \mythetag{3.3} such that 
$$
\hskip -2em
0\leqslant y_0 <2^{q_1+s_2}.
\mytag{3.10}
$$
The option $y_0=0$ is excluded since in this case we would have
$$
x_0=\frac{2^{q_1}-1}{2^{q_1+s_2}},
$$
where
$$
\hskip -2em
0<\frac{2^{q_1}-1}{2^{q_1+s_2}}<1, 
\mytag{3.11}
$$
which would mean that $x_0$ is not integer. Therefore the inequalities in \mythetag{3.10}
are rewritten as
$$
\hskip -2em
0<y_0 <2^{q_1+s_2}.
\mytag{3.12}
$$
Since $(x_0,y_0)$ is a solution of the Diophantine equation \mythetag{3.7}, we
have 
$$
\hskip -2em
x_0=\frac{3^{s_1}}{2^{q_1+s_2}}\,y_0+\frac{2^{q_1}-1}{2^{q_1+s_2}}.
\mytag{3.13}
$$
Applying \mythetag{3.12} to \mythetag{3.13}, we get
$$
\hskip -2em
\frac{2^{q_1}-1}{2^{q_1+s_2}}<x_0<3^{s_1}+\frac{2^{q_1}-1}{2^{q_1+s_2}}.
\mytag{3.14}
$$
Taking into account that $x_0$ is an integer number and taking into account
the inequalities \mythetag{3.11}, from \mythetag{3.14} we derive
$$
\hskip -2em
0<x_0\leqslant 3^{s_1}.
\mytag{3.15}
$$
Summarizing the above considerations, we can formulate the following theorem. 
\mytheorem{3.3} For any sequence of positive integers \mythetag{2.9} the Diophantine 
equation \mythetag{3.7} has a unique particular solution $(x_0,y_0)$ such that
$x_0$ and $y_0$ obey the inequalities \mythetag{3.15} and \mythetag{3.12}. 
\endproclaim
     Theorem~\mythetheorem{3.3} is similar to Item 3 in Theorem~\mythetheorem{3.1}. 
Relying on Theorem~\mythetheorem{3.3}, we can introduce the following two functions
$$
\xalignat 2
&\hskip -2em
x_0=X_0(s_1,q_1,s_2),
&&y_0=Y_0(s_1,q_1,s_2).
\mytag{3.16}
\endxalignat
$$
Substituting \mythetag{3.16} into the equation  \mythetag{3.7}, we get
$$
\hskip -2em
2^{q_1+s_2}\cdot x_0-3^{s_1}y_0=2^{q_1}-1.
\mytag{3.17}
$$
For any sequence of positive integers \mythetag{2.9} the first term in
\mythetag{3.17} is even, the coefficient $3^{s_1}$ is odd, and the right
side $2^{q_1}-1$ is also odd. Therefore $y_0$ in \mythetag{3.17} should
be odd. We have proved the following theorem.
\mytheorem{3.4} If the arguments of the function $Y_0(s_1,q_1,s_2)$ are
positive, then its value is positive and odd.
\endproclaim
     The values of the function $X_0(s_1,q_1,s_2)$ can be either even or odd.
However, solving the Diophantine equation \mythetag{3.3}, we need to get odd
numbers $m_1$ and $m_2$. Therefore we should choose proper values of $t$ in
\mythetag{3.9}:
$$
\align
&\hskip -4em
m_1=\cases Y_0(s_1,q_1,s_2)+2^{q_1+s_2}\cdot 2\,t&\text{if \ }X_0(s_1,q_1,s_2)
\text{\ \ is odd,}\\
Y_0(s_1,q_1,s_2)+2^{q_1+s_2}\cdot(2\,t+1)&\text{if \ }X_0(s_1,q_1,s_2)
\text{\ \ is even,}\\
\endcases
\mytag{3.18}\\
\vspace{1ex}
&\hskip -4em
m_2=\cases X_0(s_1,q_1,s_2)+3^{s_1}\cdot 2\,t&\text{if \ }X_0(s_1,q_1,s_2)
\text{\ \ is odd,}\\
X_0(s_1,q_1,s_2)+3^{s_1}\cdot(2\,t+1)&\text{if \ }X_0(s_1,q_1,s_2)
\text{\ \ is even,}\\
\endcases
\mytag{3.19}
\endalign
$$
The formulas \mythetag{3.18} and \mythetag{3.19} provide the required odd 
numbers $m_1$ and $m_2$ solving the equation \mythetag{3.3}. Therefore
the odd number $n=2^{s_1}\cdot m_1-1$ proves Theorem~\mythetheorem{2.2}
for the case $r=1$ in \mythetag{2.9}.\par
     The functions $X_0(s_1,q_1,s_2)$ and $Y_0(s_1,q_1,s_2)$ in \mythetag{3.16}
cannot be expressed by formulas. However they can be computed using Euclidean 
algorithm.
\medskip
\parshape 1 10pt 350pt 
{\tt\noindent
\darkred{Dioph\_solve:=proc(A,B,C) option remember:}\newline
\darkred{\ local AA,BB,XY,q:}\newline
\darkred{\ if B=0 then return [C/A,0]}\newline
\darkred{\ elif A=0 then return [0,-C/B]}\newline
\darkred{\ elif A<B}\newline
\darkred{\ \ then }\newline
\darkred{\ \ \ BB:=irem(B,A,'q'):}\newline
\darkred{\ \ \ XY:=procname(A,BB,C):}\newline
\darkred{\ \ \ return [XY[1]+q*XY[2],XY[2]]}\newline
\darkred{\ else}\newline
\darkred{\ \ AA:=irem(A,B,'q'):}\newline
\darkred{\ \ XY:=procname(AA,B,C):}\newline
\darkred{\ \ return [XY[1],XY[2]+q*XY[1]]}\newline
\darkred{\ end if}\newline
\darkred{end proc:}}
\medskip
\noindent
The above code solves the Diophantine equation $A\,x-B\,y=C$ with $A>0$ and
$B>0$. Below is the code for the function $Y_0(s_1,q_1,s_2)$.
\medskip
\parshape 1 10pt 350pt 
{\tt\noindent
\darkred{Y0:=proc(s1,q1,s2) option remember:}\newline
\darkred{\ local A,B,C,XY,YY:}\newline
\darkred{\ A:=2\^{}(q1+s2):}\newline
\darkred{\ B:=3\^{}s1:}\newline
\darkred{\ C:=2\^{}q1-1:}\newline
\darkred{\ XY:=Dioph\_solve(A,B,C):}\newline
\darkred{\ YY:=XY[2]:}\newline
\darkred{\ if YY<0}\newline
\darkred{\ \ then return YY+(iquo(abs(YY),A)+1)*A}\newline
\darkred{\ elif YY>=A}\newline
\darkred{\ \ then return YY-iquo(YY,A)*A}\newline
\darkred{\ else}\newline
\darkred{\ \ return YY}\newline
\darkred{\ end if}\newline
\darkred{end proc:}}
\medskip
\noindent
And finally we provide the code for the function $X_0(s_1,q_1,s_2)$.
\medskip
\parshape 1 10pt 350pt 
{\tt\noindent
\darkred{X0:=proc(s1,q1,s2) option remember:}\newline
\darkred{\ local A,B,C,XY,XX,YY:}\newline
\darkred{\ A:=2\^{}(q1+s2):}\newline
\darkred{\ B:=3\^{}s1:}\newline
\darkred{\ C:=2\^{}q1-1:}\newline
\darkred{\ XY:=Dioph\_solve(A,B,C):}\newline
\darkred{\ XX:=XY[1]:}\newline
\darkred{\ YY:=XY[2]:}\newline
\darkred{\ if YY<0}\newline
\darkred{\ \ then return XX+(iquo(abs(YY),A)+1)*B}\newline
\darkred{\ elif YY>=A}\newline
\darkred{\ \ then return XX-iquo(YY,A)*B}\newline
\darkred{\ else}\newline
\darkred{\ \ return XX}\newline
\darkred{\ end if}\newline
\darkred{end proc:}}
\medskip
\noindent
All of the above code is given using the programming language of the Maple 
package, version 9.01. Maple is a trademark of Waterloo Maple Inc.\par
\head
4. Long evolution sequences.
\endhead
     Let's proceed to the case $r=2$. In this case we have the long evolution
sequence $s_1,\,q_1,\,s_2,\,q_2,\,s_3$ in \mythetag{2.9} that subdivides into two 
short sequences
$$
\xalignat 2
&\hskip -2em
s_1,\,q_1,\,s_2, 
&&s_2,\,q_2,\,s_3.
\mytag{4.1}
\endxalignat
$$
The sequences \mythetag{4.1} generate two Diophantine equations similar to 
\mythetag{3.7}. The solution of the first one is given by the formulas
\mythetag{3.18} and \mythetag{3.19}. Let's write these two formulas as 
follows:
$$
\xalignat 2
&\hskip -2em
m_1=M_{10}+2^{q_1+s_2}\cdot 2\,t_{1},
&&m_2=\tilde M_{20}+3^{s_1}\cdot 2\,t_{1}.
\mytag{4.2}
\endxalignat
$$
The solution of the second Diophantine equation is given by similar
formulas
$$
\xalignat 2
&\hskip -2em
m_2=M_{20}+2^{q_2+s_3}\cdot 2\,t,
&&m_3=\tilde M_{30}+3^{s_2}\cdot 2\,t.
\mytag{4.3}
\endxalignat
$$
The formulas \mythetag{4.2} and \mythetag{4.3} produce a new
Diophantine equation
$$
\hskip -2em
2^{q_2+s_3}\cdot t-3^{s_1}\cdot t_1=\frac{\tilde M_{20}-M_{20}}{2}
\mytag{4.4}
$$
with respect to the variables $t_1$ and $t$. Note that $\tilde M_{20}$ 
and $M_{20}$ are odd. Therefore the value of the fraction in the right 
hand side of the equation \mythetag{4.4} is integer.\par 
    The equation \mythetag{4.4} is very similar to \mythetag{3.7}. It is always 
solvable due to Item 1 of Theorem~\mythetheorem{3.1} since $\GCD(2^{q_2+s_3},3^{s_1})=1$. 
Its solution is written as
$$
\xalignat 2
&\hskip -2em
t_{1}=T_{20}+2^{q_2+s_3}\cdot t_2,
&&t=\tilde T_{21}+3^{s_1}\cdot t_2.
\mytag{4.5}
\endxalignat
$$
Substituting \mythetag{4.5} into the first equality \mythetag{4.2} and 
into the second equality \mythetag{4.3}, we derive
$$
\xalignat 2
&\hskip -2em
m_1=M_{11}+2^{q_1+q_2+s_2+s_3}\cdot 2\,t_2,
&&m_3=\tilde M_{31}+3^{s_1+s_2}\cdot 2\,t_2,
\mytag{4.6}
\endxalignat
$$
\vskip -1ex
\noindent
where
\vskip -5ex
$$
\xalignat 2
&\hskip -2em
M_{11}=M_{10}+2^{q_1+s_2}\cdot 2\,T_{20},
&&\tilde M_{31}=\tilde M_{30}+3^{s_2}\cdot 2\,\tilde T_{21}.
\mytag{4.7}
\endxalignat
$$
The formulas \mythetag{4.6} and \mythetag{4.7} serve the case $r=2$ in
\mythetag{2.9}. The odd number $n=2^{s_1}\cdot m_1-1$ proves 
Theorem~\mythetheorem{2.2} for this case.\par 
     The formulas \mythetag{4.6} are similar to \mythetag{4.2}. Therefore we can
increment $r$ by $1$, complement the sequences \mythetag{4.1} by $s_3,\,q_3,\,s_4$, 
and write the formulas
$$
\xalignat 2
&\hskip -2em
m_3=M_{30}+2^{q_3+s_4}\cdot 2\,t,
&&m_4=\tilde M_{40}+3^{s_3}\cdot 2\,t.
\mytag{4.8}
\endxalignat
$$
The equalities \mythetag{4.6} and \mythetag{4.8} produce a new Diophantine equation
$$
\hskip -2em
2^{q_3+s_4}\cdot t-3^{s_1+s_2}\cdot t_2=\frac{\tilde M_{31}-M_{30}}{2}
\mytag{4.9}
$$
similar to \mythetag{4.4}. The equation \mythetag{4.9} is similar to \mythetag{3.7}. 
It is always solvable due to Item 1 of Theorem~\mythetheorem{3.1} since 
$\GCD(2^{q_3+s_4},3^{s_1+s_2})=1$. Its solution is written as
$$
\xalignat 2
&\hskip -2em
t_2=T_{30}+2^{q_3+s_4}\cdot t_3,
&&t=\tilde T_{31}+3^{s_1+s_2}\cdot t_3.
\mytag{4.10}
\endxalignat
$$
We substitute \mythetag{4.10} into the first equality \mythetag{4.6} and 
into the second equality \mythetag{4.8}. As a result we get the equalities 
$$
\xalignat 2
&\hskip -2em
m_1=M_{12}+2^{q_1+q_2+q_3+s_2+s_3+s_4}\cdot 2\,t_3,
&&m_4=\tilde M_{41}+3^{s_1+s_2+s_3}\cdot 2\,t_3,
\qquad
\mytag{4.11}
\endxalignat
$$
\vskip -1ex
\noindent
where
\vskip -5ex
$$
\xalignat 2
&\hskip -2em
M_{12}=M_{11}+2^{q_1+q_2+s_2+s_3}\cdot 2\,T_{30},
&&\tilde M_{41}=\tilde M_{40}+3^{s_3}\cdot 2\,\tilde T_{31}.
\mytag{4.12}
\endxalignat
$$\par
     Continuing the process we can go further by induction. 
Extending the sequence of formulas \mythetag{4.2}, \mythetag{4.6}, 
and \mythetag{4.11}, we write
$$
\xalignat 2
&\hskip -2em
m_1=M_{1\,r-1}+2^{\,Q_r}\cdot 2\,t_r,
&&m_{r+1}=\tilde M_{r+1\,1}+3^{S_r}\cdot 2\,t_r
\mytag{4.13}
\endxalignat
$$
The exponentials $Q_r$ and $S_r$ are given by the formulas
$$
\xalignat 2
&\hskip -2em
Q_r=\sum^r_{i=1}q_i+\sum^{r+1}_{i=2}s_i,
&&S_r=\sum^r_{i=1}s_i.
\mytag{4.14}
\endxalignat
$$
They obey the relationships which are used below in proving the inductive step:
$$
\xalignat 2
&\hskip -2em
Q_r=Q_{r-1}+q_r+s_{r+1},
&&S_r=S_{r-1}+s_r.
\mytag{4.15}
\endxalignat
$$
The recurrent relationships \mythetag{4.15} are immediate from 
\mythetag{4.14}.\par
     For $r>1$ the quantities $M_{1\,r-1}$ and $\tilde M_{r+1\,1}$ in 
\mythetag{4.13} are defined inductively:
$$
\align
&\hskip -2em
M_{1\,r-1}=M_{1\,r-2}+2^{\,Q_{r-1}}\cdot 2\,T_{r0}.
\mytag{4.16}\\
&\hskip -2em
\tilde M_{r+1\,1}=\tilde M_{r+1\,0}+3^{s_r}\cdot 2\,\tilde T_{r1},
\mytag{4.17}
\endalign
$$
where $T_{r0}$ and $\tilde T_{r1}$ are defined by solving the Diophantine equation
$$
\hskip -2em
2^{q_r+s_{r+1}}\cdot t-3^{S_{r-1}}\cdot t_{r-1}=\frac{\tilde M_{r1}-M_{r0}}{2}
\mytag{4.18}
$$
whose general solution is taken in the following form:
$$
\xalignat 2
&\hskip -2em
t_{r-1}=T_{r0}+2^{q_r+s_{r+1}}\cdot t_r,
&&t=\tilde T_{r1}+3^{S_{r-1}}\cdot t_r.
\mytag{4.19}
\endxalignat
$$\par
      In \mythetag{4.18} and in \mythetag{4.17} we see the quantities 
$M_{r0}$ and $\tilde M_{r+1\,0}$ respectively, \mythetag{4.12} being 
a particular case of \mythetag{4.18} and \mythetag{4.17} for $r=3$. These 
quantities are determined by the short sequence of positive integers 
$s_r,\,q_r,\,s_{r+1}$ through the functions $X_0$ and $Y_0$ defined in the 
previous section: 
$$
\align
&\hskip -4em
M_{r0}=\cases Y_0(s_r,q_r,s_{r+1})&\text{if \ }X_0(s_r,q_r,s_{r+1})
\text{\ \ is odd,}\\
Y_0(s_r,q_r,s_{r+1})+2^{q_r+s_{r+1}}&\text{if \ }X_0(s_r,q_r,s_{r+1})
\text{\ \ is even,}\\
\endcases
\mytag{4.20}\\
\vspace{1ex}
&\hskip -4em
\tilde M_{r+1\,0}=\cases X_0(s_r,q_r,s_{r+1})&\text{if \ }X_0(s_r,q_r,s_{r+1})
\text{\ \ is odd,}\\
X_0(s_r,q_r,s_{r+1})+3^{s_r}&\text{if \ }X_0(s_r,q_r,s_{r+1})
\text{\ \ is even.}\\
\endcases
\mytag{4.21}
\endalign
$$
The formulas \mythetag{4.20} and \mythetag{4.21} mean that $M_{r0}$ and 
$\tilde M_{r+1\,0}$ are always odd.\par
     The quantities $M_{1\,r-1}$ and $\tilde M_{r+1\,1}$ in \mythetag{4.13} are 
defined inductively by means of the formulas \mythetag{4.16} and \mythetag{4.17}
for $r>1$. The case $r=1$ is the base of this induction. In this case $M_{1\,r-1}$
turns to $M_{10}$ which is given by the formula \mythetag{4.20}. The case $r=1$
is the base of induction for the formula \mythetag{4.13} as well. In this case
$\tilde M_{r+1\,1}$ turns to $\tilde M_{21}$, while \mythetag{4.13} turns
to \mythetag{4.2}. Therefore $\tilde M_{21}=\tilde M_{20}$ and $\tilde M_{20}$
is given by the formula \mythetag{4.21}.\par
     In order to prove the formulas \mythetag{4.13} now it is sufficient to prove the
inductive step $r-1\to r$. Let's replace $r$ by $r-1$ in \mythetag{4.13} and assume 
that the formulas obtained from \mythetag{4.13} in such a way are valid:
$$
\xalignat 2
&\hskip -2em
m_1=M_{1\,r-2}+2^{\,Q_{r-1}}\cdot 2\,t_{r-1},
&&m_{r}=\tilde M_{r\,1}+3^{S_{r-1}}\cdot 2\,t_{r-1}
\mytag{4.22}
\endxalignat
$$
We complement these formulas with the formulas similar to \mythetag{4.2} and
\mythetag{4.3}: 
$$
\xalignat 2
&\hskip -2em
m_r=M_{r0}+2^{q_r+s_{r+1}}\cdot 2\,t,
&&m_{r+1}=\tilde M_{r+1\,0}+3^{s_r}\cdot 2\,t.
\mytag{4.23}
\endxalignat
$$
These formulas are derived from the Diophantine equation similar to \mythetag{3.7} 
and associated with the short sequence $s_r,q_r,s_{r+1}$.\par
    The second formula \mythetag{4.22} and the first formula \mythetag{4.23} 
represent the same quantity $m_r$. Equating their right hand sides, we obtain a 
Diophantine equation with respect to $t_{r-1}$ and $t$. This Diophantine equation
coincides with \mythetag{4.18}. Its solution is given by the formulas \mythetag{4.19}.
The formulas \mythetag{4.13} then are derived by substituting \mythetag{4.19} into
the first formula \mythetag{4.22} and into the second formula \mythetag{4.23} if
we take into account \mythetag{4.15}, \mythetag{4.16}, and \mythetag{4.17}.\par
     Thus, the formulas \mythetag{4.13} are proved. They serve the general case
$r>1$ in \mythetag{2.9}. The first formula \mythetag{4.13} determines the positive
odd number $m_1$ depending on an arbitrary positive integer parameter $t_r$. 
For any integer value of this parameter the positive odd number $n=2^{s_1}\cdot m_1-1$
proves Theorem~\mythetheorem{2.2} in the case of an arbitrary finite sequence
of positive integers \mythetag{2.9}.\par 
\head
5. Conclusions.
\endhead
     Theorem~\mythetheorem{2.2} is the main result of the present paper. It positively 
solves the problem of a predefined Collatz evolution for finite length evolution patterns 
\mythetag{2.9}. But this result cannot be easily transferred to the case of infinite
patterns \mythetag{2.7}.
\head
6. Dedicatory.
\endhead
     This paper is dedicated to my sister Svetlana Abdulovna Sharipova. 

\enddocument
\end